\documentclass{amsart}

\usepackage{amsthm,amsfonts,amsmath,amssymb,mathrsfs,verbatim,cite}

\numberwithin{equation}{section}
\newtheorem{theorem}{Theorem}[section]
\newtheorem{corollary}[theorem]{Corollary}
\newtheorem{proposition}[theorem]{Proposition}

\newtheorem{lemma}[theorem]{Lemma}

\newtheorem{definition}[theorem]{Definition}
\newtheorem{example}[theorem]{Example}

\newcommand{\Real}{\mathbb R}

\newcommand{\Nat}{\mathbb N}
\newcommand{\Z}{\mathbb Z}

\DeclareMathOperator{\Nil}{Nil}

\begin{document}

\title{Right $P$-comparable semigroups}

\author{Nazer. H. Halimi}

\thanks{}

\email{n.halimi@uq.edu.au}

\address{School of Mathematics and Physics, 
The University of Queensland, QLD 4072, Australia}
\date{November 2009}

\subjclass[2000]{Primary 20M12; Secondary 12M11}

\begin{abstract}
In this paper we introduce the notion of right waist and right comparizer ideals for semigroups. 
In particular, we study the ideal theory of semigroups containing right waists and right comparizer ideals. 
We also study those properties of right cones that can be carried over to right 
$P$-comparable semigroups. We give sufficient and necessary conditions 
on the set of nilpotent elements of a semigroup to be an ideal. 
We provide several equivalent characterizations for a right ideal being 
a right waist. In one of our main result we show that in a right $P_1$-comparable semigroup 
with left cancellation law, a prime segment $P_2\subset P_1$ is 
Archimedean, simple or exceptional, extending a similar result of right cones to $P$-comparable semigroups. 
\end{abstract}

\maketitle 
\section{Introduction} 

In ring and module theory the notion of a \textit{waist} was first
introduced by Auslander \textit{et al.} to study a new class of 
indecomposable modules, which they called modules with waists \cite{AGR75}. 
In particular, new modules were constructed from known modules with emphasis on 
modules with waists.
In \cite{FT93}, Ferrero and Torner presented a complete 
characterization of right waists contained in the Jacobson radical of 
a ring $R$. In \cite{FT95} these same authors studied prime $D$-rings 
with \textit{(MP)} property satisfying an ascending chain condition on right waists,
and provided a complete characterization of right $D$-domains of such type.
Mazurek and Torner in \cite{MT06} used \textit{right comparizers} 
and right waists to provide a characterization theorem for semiprime 
segments of a ring. 
In 1999  \textit{right $P$-comparable} rings were 
introduced by Ferrero and Sant' Ana \cite{FSA99}. Here $P$ is a completely prime right ideal,
 which are proved in \cite{FSA99} to be right waists. Examples of right $P$-comparable rings
are right distributive rings with \textit{(MP)} property.
These authors studied prime and semiprime ideals, right Noetherian rings with 
comparability and provided a structure theorem for such rings.   

The purpose of this paper is to extend the above-mentioned developments
to the realm of semigroup theory. In particular we study the ideal theory 
of semigroups containing analogous notions of
right waist and right comparizer ideals, and investigate those properties 
of right cones that can be carried over to right $P$-comparable semigroups. 
Among our results we obtain a necessary and sufficient condition on the set 
of nilpotent elements of a semigroup to be an ideal. Furthermore, we prove
that right $P$-comparability is equivalent to weak right $P$-comparability. 
We also show that in a right $P_1$-comparable semigroup with left 
cancellation law, a prime segment $P_2\subset P_1$ is Archimedean, 
simple or exceptional, extending the analogous result from right cones.
\medskip
  
In the first part of the paper we investigate the properties of right waist 
and right comparizer ideals in semigroups and transfer known results and 
ideas from ring theory to semigroups theory. 
The papers \cite{FT93} and \cite{MT04} were essential for the this  
part of the paper. In the second part of this work we generalize 
properties of ideal theory in right cones to right $P$-comparable 
semigroups with left cancellation law.

\medskip

Section~\ref{Sec 2}, concerns definitions and properties of right waist 
and right comparizer ideals and their relation with classical radicals. 
For example; in Theorem~\ref{Thm 2.10}, by using right comparizer 
ideals, we give a characterization theorem for the set of nilpotent elements 
to form an ideal.

In sections~\ref{Sec 3} and \ref{Sec 4} we focus on right 
$P$-comparable semigroups. In the first of these sections we investigate the 
properties of prime and completely prime ideals of right 
$P$-comparable semigroups. In Theorem~\ref{Thm 3.8}, it is shown that right 
$P$-comparability is equivalent to right weak $P$-comparability. 
In section~\ref{Sec 4}, prime segments of right $P$-comparable 
semigroups are investigated. We show that right $P$-comparable 
semigroups have many properties in common with right cones. 
One of our main result is that in a right $P$-comparable 
semigroup $S$ with left cancellation law, a prime segment is either 
Arthimedean, simple or exceptional. 

\section{Right Waist and right Comparizer ideals of semigroups}\label{Sec 2}
Throughout this paper, $S$ always denotes a semigroup with unit element $1$ and zero element $0$ such that $1\neq 0$. 
These semigroups are not necessary commutative.  

A proper right ideal $P$ of a semigroup $S$ is called \textit{prime}, 
if for any $a,b\in S$, $aSb\subseteq P$ implies either $a\in P$ or $b\in P$. 
If $ab\in P$ implies either $a\in P$ or $b\in P$, then $P$ is called \textit{completely prime}.
Moreover, $P$ is said to be \textit{semiprime}, if $aSa\subseteq P$ implies $a\in P$.

In analogy with ring theory, a proper right ideal $I$ of a semigroup $S$ 
is said to be a \textit{right waist}, if $I$ is comparable with every right 
ideal of $S$, that is, either $A\subseteq I$ or $I\subset A$ holds for any 
right ideal $A$ of $S$. A right ideal $I$ of a semigroup $S$ is a 
\textit{right comparizer} if for all right ideals $A,B$ of $S$, either 
$A\subseteq B$ or $BI\subseteq A$. This is equivalent to the statement that, 
for any $a,b\in S$, either $aS\subseteq bS$ or $bI\subseteq aS$.

A semigroup $S$ is \textit{right $P$-comparable} with respect to a 
completely prime right ideal $P$, if for every $a, b\in S$ one of 
the following conditions hold; 
$aS\subseteq bS$, $bS\subseteq aS$ or $(aS)T^{-1}=(bS)T^{-1}$, where $T=S-P$. 

From the definition of right comparizer ideal we immediately obtain the following.
\begin{lemma}\label{Le 2.1}

For a semigroup $S$, the following right ideals are right comparizer:

\begin{itemize}
\item[(i)] The zero ideal of $S$.

\item[(ii)] The union of a family of right comparizer ideals of $S$.

\item[(iii)] A right ideal of $S$ contained in a right comparizer ideal of $S$. 

\end{itemize}
\end{lemma}

The following lemma shows a similarity between right waist idempotent and right waist completely prime ideals.
\begin {lemma}\label{Le 2.2}
Let $S$ be a semigroup. Then the following statements hold:

\begin{itemize}  

\item[(i)] If $I$ is an idempotent ideal of $S$ and a right waist, then $I=aI$ for every $a\in S-I$.

\item[(ii)] A completely prime ideal $P$ is a right waist if and only if $P=aP$ for every $a\in S-P$.
\end{itemize}
\end {lemma}

\begin {proof}
(i) Since $I$ is a right waist, we have $I\subseteq aS$ for every $a\in S-I$. Therefore, $I=I^{2}\subseteq aSI=aI$, which shows that $I\subseteq aI$.
 The converse is clear so that $I=aI$.

(ii) Let $I$ be a right ideal and $a\in I-P$. Then $P=aP\subseteq I$, and $P$ is a right waist. Conversely, since $P$ is a right waist,
 $P\subseteq aS$ for every $a\in S-P$. Therefore, for each $p\in P$ there exists $s\in S$ with $p=as$. Since $a\notin P$, $s\in P$ and $P\subseteq aP\subseteq P$.
\end {proof}

An element $u\in S$ is called a unit if there exists an element $v\in S$ such that $uv=vu=1$,\; $U=U(S)$ stands for the group of units.
The set $J(S)=S-U(S)$ is  called the \textit{non units} of $S$. 
For any $a,b,c\in S$, if $ab=ac\neq 0$ implies $b=c$, then we say that $S$ has a left cancellation law. 
Right cancellation law is defined similarly. A semigroup $S$ has cancellation law provided $S$ has left and right cancellation laws.

\begin{proposition}\label{Pr 2.3}
Let $S$ be a cancellation semigroup. Then:

\begin{itemize}
\item[(i)]
The non units $J(S)$ is a maximal right and maximal left ideal of $S$. 

\item[(ii)]
 The non units $J(S)$ is a completely prime ideal.

\item[(iii)]

 $J(S)=\cup I$, where $I$ are proper right ideals of $S$.

\item[(iv)]
$J(S)=\cup I$, where $I$ are proper left ideals of $S$.
\end{itemize}
\end{proposition}
\begin{proof}
(i)
Suppose that $j\in J(S), s\in S$ and $js\notin J(S)$. Then $js\in U(S)$ and $(js)u=j(su)=1$ for some $u\in U(S)$. Also $j=1j=jsuj$. 
Hence $1=(su)j$, a contradiction. Similarly it can be shown that $J(S)$ is left ideal of $S$.
 If $M$ is a left ideal of $S$ with $J(S)\subset M$, then there exists an element $m\in M-J(S)$ which is a unit. 
This shows that $M=S$. The same reasoning shows that $J(S)$ is a maximal right ideal of $S$.

(ii)
 If $a,b\notin J(S)$, then $a,b\in U(S)$ and $ab\in U(S)$. Hence $J(S)$ is a completely prime ideal.

(iii)
 Since $J(S)$ is a proper right ideal of $S$, $J(S)\subseteq \cup I$. Conversely, if $a\in \cup I$, then $a\in I\neq S$ and $a\notin U(S)$.

The proof of (iv) is similar to (iii).
\end{proof}
 
\begin{theorem}\label{Thm 2.4}

Let $I$ be a right comparizer ideal of a semigroup $S$. 
 
 Then the following three statements hold:
 \begin{itemize}

\item[(i)] If $I=I^{2}$, then $I$ is a right waist.

\item[(ii)] If $I$ is a right waist, then $aI$ is a right waist for every $a\in S$.

\item[(iii)] If $P$ is a prime right ideal of $S$ such that $I\nsubseteq P$, then $P$ is a right waist and $P\subset I$.
 Hence the set of prime right ideals contained in $I$ is linearly ordered by inclusion.
\end{itemize}
Furthermore, if $S$ has left cancellation law, then the following two statements hold:
\begin{itemize}
\item[(iv)]
 If $I$ is a nonnilpotent ideal, then $\cap_{n\in \Nat} I^{n}$ is completely prime.

\item[(v)]
 If $I$ is an idempotent ideal, then $I$ is completely prime.
 \end{itemize}
\end{theorem}

\begin {proof}
(i) Let $A$ be a right ideal of $S$. Since $I$ is right comparizer, we have $A\subseteq I$ or $I=I^{2}\subseteq A$. 
Hence $I$ is a right waist. To prove (ii), let $b\in S$ and $bS\nsubseteq aI$. We need to show that $aI\subseteq bS$. 
Since $I$ is right comparizer, we have $bS\subseteq aS$ or $aI\subset bS$. 
Thus it is enough to consider $bS\subseteq aS$, which shows that $b=as$ for some $s\in S$.
 By assumption $b\notin aI$, so that $s\notin I$. Since $I$ is a right waist, $I\subseteq sS$, and $aI\subseteq asS=bS$.

(iii) Let $A$ be a right ideal of $S$ such that $P\nsubseteq A$. Then $AI\subseteq P$ and $A\subseteq P$ follows. Thus $P$ is a right waist.

(iv)
If $I=S$ then there is nothing to prove. We can thus assume that $I$ is a 
proper ideal of $S$ and $I\subseteq J(S)$.  
We proceed contradiction, and assume that $ab\in \cap_{n\in\Nat}I^n$
for some $a,b\in S-\cap_{n\in\Nat}I^n$. Then $a,b\notin I^{m}$ for 
some positive integer $m$.  
Since $aS\nsubseteq I^m$, $bS\nsubseteq I^m$ and $I$ is right comparizer, 
we have $I^{m+1}\subseteq aS$ and $I^{m+1}\subseteq bS$. 
Thus $I^{2m+2}\subseteq  aSI\cap bSI\subseteq aI\cap bI$ and so
$I^{2m+3}=I^{m+2}I^{m+1}\subseteq aII^{m+1}=aI^{m+2}\subseteq 
abI\subseteq abJ(S)$. 
Since $ab\in\cap I^n\subseteq I^{2m+3}$, we obtain $ab\in abJ(S)$ and,  
by the left cancellation law, we have $1\in J(S)$, a contradiction.

(v)
Since $I=I^2=I^n$ for all $n$, we have $I=\cap I^{n}$. Hence by (iv), 
$I$ is completely prime.
\end {proof}

The following lemma provides some examples of right comparizer ideals. 
\begin{lemma}\label{Le 2.5}

Let $I$ be a right waist of a semigroup $S$. Then: 
\begin{itemize}

\item[(i)] A right ideal $C$ of $S$ contained in $I$ is a comparizer right ideal of $S$ if and only if $C$ is comparizer among all right ideals $S$ contained in $I$.

\item[(ii)] $I\cap r_{S}(I)$ is a comparizer ideal of $S$, where $r_{S}(I)$ is the right annihilator of $I$.

\item[(iii)] If $I^{n}=0$ for $n\geq 2$ then $I^{n-1}$ is a comparizer ideal of $S$.

\end{itemize}
\end{lemma}

\begin{proof} 
The proof is the same as that in ring theory, see \cite[Lemma 1.3]{MT04}.
\end{proof}

For a semigroup $S$, let $C(S)$ denote the union of all proper right comparizer ideals of $S$. Then $C(S)$ is the largest right comparizer ideal of $S$. 
Moreover every right ideal of $S$ contained in $C(S)$ is right comparizer. From the definition we obtain the following lemma.

\begin{lemma}\label{Le 2.6}
Let $S$ be a semigroup. Then: 
\begin{itemize}
\item[(i)]
$C(S)=\{c\in S: \text{for all $a,b\in S$, either $a\in bS$ or $bc\in aS$}\}$.

\item[(ii)]
$S$ is a right chain semigroup if and only if $C(S)=S$.
\end{itemize}
\end{lemma}

Next we investigate the relation between right comparizer and right waist ideals and the following five \textit{radicals}. 

\begin{itemize}

\item
The prime radical $\beta(S)$ define as the intersection of all prime ideal of $S$.

\item
$N(S)$, the intersection of all completely prime ideals of $S$.

\item
The nil radical ${\Nil}(S)$ define as the largest nil ideal of $S$. 

\item
$A(S)$, the union of all nilpotent ideals of $S$.

\item
$T(S)$, the set of all nilpotent elements of $S$.

\end{itemize}

\begin{theorem}\label{Thm 2.7}
Let $S$ be a semigroup. Then the following holds:
\begin{itemize}

\item[(i)]
 If $C(S)$ is nilpotent, then $C(S)\subseteq N(S)$.

\item[(ii)]
If $C(S)$ is nonnilpotent, then $N(S)\subseteq C(S)$ and $N(S)$ is a completely prime ideal and a right waist.

\item[(iii)]
If $C(S)$ is nonnilpotent, then $\beta(S)$ is prime and a right waist.
\end{itemize}
\end{theorem}
\begin{proof}

The proof of (i) is obvious. To prove (ii) we first assume that $P\subseteq C(S)$ for some completely prime ideals of $S$.  
Then, by part (iii) of the Theorem~\ref{Thm 2.4}, $N(S)$ is a completely prime ideals and a right waist.
 Next we assume that  $P\nsubseteq C(S)$ for all completely prime. Since $C(S)$ is right comparizer, 
we have $C(S)^2\subseteq P$, and so $C(S)\subset P$ for all completely ideals of $P$. 
By the fact that $C(S)$ is nonnilpotent and using of part (iv) of Theorem~\ref{Thm 2.4}, we have $N(S)\subseteq \cap_nC(S)^n\subseteq C(S)\subseteq N(S)$.
 Therefore, $N(S)=\cap_nC(S)^n$ is completely prime and also a right waist.
The case (iii) is obvious from (ii) and part (iii) of Theorem~\ref{Thm 2.4}.
\end{proof}

The following theorem is an extension of Theorem 6 of Ferrero \textit{et al.} \cite{FMSA05} to semigroups with left cancellation law.
\begin{theorem}\label{Thm 2.8}

Let $S$ be a left cancellation semigroup with nonnilpotent $C(S)$. Then:
\begin{itemize} 

\item[(i)]
If $t\in T(S)$ and $a\in S$, then $taP\subseteq aP$.

\item[(ii)]
$T(S)$ is a subsemigroup (without identity) of $S$, and $T(S)$ is equal to the union of the nilpotent ideals of $T(S)$.

\item[(iii)]
$A(S)=\beta(S)=\Nil(S)$ is a prime ideal and a right waist.
\end{itemize}
\end{theorem}
\begin{proof}
To prove (i) 
we assume that $taP\nsubseteq aP$. Since $aP$ is a right waist, we have $aP\subseteq taP$ and $aP\subseteq t^{n}aP$ for all positive integer $n$.
 Since $t\in T(S)$, we obtain $aP=0$, a contradiction.
 
(ii)
We put $P=N(S)$. Obviously, $T(S)\subseteq P\subseteq C(S)$, $P$ 
is completely prime and a right waist.  
If $P$ is nilpotent, then $T(S)=P\subseteq J(S)$ and we are done. 
We can thus assume that $P$ is nonnilpotent, and so $P^2=P$ by part (iv) of Theorem~\ref{Thm 2.4}.  
Let $t, t_1,t_2,\dots,t_n\in T(S)$ with $t^n=0$. 
From (i) one can easily show that $tt_1 tt_2 \cdots tt_n 
P\subseteq t^nP=0$. Hence $(tT(S))^n=0$ as desired.

(iii)
 Let $t\in{\Nil}(S)$. Then $(StS)^3\subseteq {\Nil}(S)t{\Nil}(S)\subseteq 
T(S)tT(S)$ and from (ii) we deduce that $StS$ is a nilpotent ideal of $S$. 
Hence $t\in A(S)$ and thus $A(S)=\beta(S)={\Nil}(S)$. 
By part (iii) of Theorem~\ref{Thm 2.4}, $\beta(S)$ is a prime ideal and a right waist. 
\end{proof}

\begin{corollary}\label{Co 2.9}
  
If $S$ is a left cancellation semigroup with nonnilpotent $C(S)$, then for every ideal $I$ of $S$, either $I\subseteq\beta(S)$ or $N(S)\subseteq I$.
\end{corollary}
\begin{proof}
If $I$ is an ideal of $S$ such that $\beta(S)\subset I\subset N(S)$, then $I$ is nonnilpotent.
 By part (iv) of Theorem~\ref{Thm 2.4}, $N(S)$ is contained in $I$, a contradiction.
\end{proof}

In 1970 Skornyakov posed the question whether, in a left cone $S$
with cancellation law such that $U(S)={1}$, the set of nilpotent elements forms an ideal.
In \cite{BT98} Brungs and Torner answered Skornykov's question and characterized that in a right cone
the set of nilpotent elements forms an ideal if and only if $\beta(S)$ is completely prime.
In \cite{FMSA05} Ferrero, Mazurek and Sant'~Ana generalized this result to right chain semigroups. 
Here we prove that this result also holds for left cancellation semigroups $S$ with nonnilpotent $C(S)$.
\begin{theorem}\label{Thm 2.10}

Let $S$ be a left cancellation semigroup with nonnilpotent $C(S)$. Then the following conditions are equivalent.
\begin{itemize}

\item[(i)]
$T(S)$ is an ideal of $S$.

\item[(ii)]
 $T(S)=\beta(S)$.

\item[(iii)]
$\beta(S)$ is a completely prime ideal of $S$.
\end{itemize}
\end{theorem}
\begin{proof}

If $T(S)$ is an ideal of $S$, then $T(S)={\Nil}(S)$, and (i) implies (ii) by virtue of Theorem~\ref{Thm 2.8} part (iii). 
Moreover, $\beta(S)$ is prime. We proceed by contradiction, and assume that $\beta(S)=Q$ is an exceptional ideal; i.e., it is prime but not completely prime. 
Then $Q$ is nilpotent by part (iv) of Theorem~\ref{Thm 2.4} and part (ii) of Theorem~\ref{Thm 2.7}. 
Let $P$ be the minimal completely prime ideal of $S$ containing $Q$. Since $C(S)$ in nonnilpotent such $P$ exists, by part (iv) of Theorem~\ref{Thm 2.4}. 
As $Q$ is a right waist, and $P$ is minimal completely prime over $Q$, again by using part (iv) of Theorem~\ref{Thm 2.4},
 for every $a\in P-Q$ there exits $n$ such that $a^n\in Q$. Therefore, $a$ is a nilpotent element of $S$ and $P=T(S)$, a contradiction.
 Hence, $\beta(S)$ is a completely prime ideal of $S$, and (ii)$\Rightarrow$(iii) follows. 
If $\beta(S)$ is completely prime, then every nilpotent element of $S$ belongs to $\beta(S)$, and $\beta(S)=T(S)$. Thus $T(S)$ is an ideal of $S$.
\end{proof}

\begin{definition}\label{Def 2.11}
Let $A$ be a proper right ideal of a semigroup $S$. 
Then the set 
\[
P_r(A)=\{s\in S:~xs\in A \text{ for some } x\in S-A\}
\]
is called the \emph{associated prime right ideal} of $A$.
\end{definition}

\begin{lemma}\label{Le 2.12}

Let $A$ be a proper right ideal of a right semigroup $S$. Then:
\begin{itemize}

\item[(i)]
 $P_{r}(A)$ is a completely prime right ideal of $S$.

\item[(ii)]
 If $A$ is a prime right ideal of $S$ then for any right waist 
$I$ of $S$ either $I\subseteq A$ or $P_r(A)\subseteq I$.
\end{itemize}
\end{lemma}
\begin{proof}

(i) It is obvious from the definition that $P_r(S)$ is a right ideal of $S$. If $ab\in P_{r}(A)$, then $sab\in A$ for some $s\in S-A$.
 If $sa\in S-A$, then $b\in P_{r}(A)$, otherwise $a\in P_{r}(A)$.

(ii) Assume $P_r(A)\nsubseteq I$ and let $s\in P_r(A)-I$. Then there exists $x\in S-A$ with $xs\in A$, and $I\subseteq sS$. 
Hence $xI\subseteq xsS\subseteq A$, since $x\notin A$ and  $A$ is prime $I\subseteq A$. 
\end{proof}

Analogous of the following lemma and corollary were proved by Ferrero and Torner for $D$-rings 
\cite[Lemma 2.3 \& Corollary 2.4]{FT93}. 
\begin{lemma}\label{Le 2.13}
Let $S$ be a semigroup, $T$ be a right ideal of $S$ and $P=P_r(T)$. Then the following condition are equivalent:
\begin{itemize}
\item[(i)] 
$T$ is a right waist.

\item[(ii)]
 $T=\cap_{a\notin T}aP$ or there exists $b\in S-T$ such that $T\subset bS=\cap_{a\notin T}aP$ and $T$ is a lower neighbour of $bS$.
 \end{itemize}
\end{lemma}
\begin{proof}
To prove the implication (i)$\Rightarrow$(ii), take $a\notin T$. Since $T$ is a right waist and $aP$ 
is a right ideal we have $T\subseteq aP$ or $aP\subset T$. 
Assume there exists $b\in T-aP$. Then $b\in T\subset aS$ and 
so $b=as$ for some $s\in S$. Since $b\in T$ we get $s\in P=P_r(T)$, 
which is a contradiction because $b\notin aP$. 
Therefore $T\subseteq \cap_{a\notin T}aP$. Now if 
$b\in (\cap_{a\notin T}aP)-T$, 
then $T\subset bS\subseteq \cap_{a\notin T}aP\subseteq bP\subseteq bS$, 
so $T\subset bS=\cap_{a\notin T}aP$. 
Let $H$ be a right ideal with $T\subset H\subseteq bS$ and $h\in H-T$. 
As above we get $hS=\cap_{a\notin T}aP$ and so $H=bS$.

To obtain (i) from (ii), assume that $T\subseteq \cap_{a\notin T}aP$ 
and $H$ is a right ideal $S$ with $H\nsubseteq T$. 
Let $h\in H-T$. Then $T\subseteq hP\subseteq hS\subset H$ and 
so $T$ is a right waist.

\end{proof}
\begin{corollary}\label{Co 2.14}
Let $T$ be a right waist with $P=P_r(T)\subseteq J=J(S)$. 
Then $T=\cap_{a\notin T}aP=\cap_{a\notin T}aJ$.

\end{corollary}
\begin{proof}
If $b\in (\cap_{a\notin T}aP)-T$, then $bS=\cap_{a\notin T}aP$, and $bS\subseteq bP\subseteq bJ\subset bS$. 
This implies that $b\in bJ$, a contradiction. Thus $T=\cap_{a\notin T}aP$. If $b\in \cap_{a\notin T}aJ-T$, then we again reach the contradiction $b\in bJ$.
 Since $T$ is a right waist, $\cap_{a\notin T}aJ\subseteq T=\cap_{a\notin T}aP\subseteq \cap_{a\notin T}aJ$.
\end{proof}

\section{right $P$-comparable semigroups}\label{Sec 3}

Let $S$ be a semigroup. For every multiplicatively closed subset $T$ of $S$ and $a\in S$, we can define the set $(aS)T^{-1}:=\{x\in S$: $xt\in (aS)$ for some $t\in T \}$.
 Recall that $T$ is said to be a right \textit{Ore} set if for every $a\in S$ and $t\in T$ there exist $a'\in S$ and $t'\in T$, such that $at'=ta'$.

\begin{lemma}\label{Le 3.1}
Assume that $T\subseteq S$ is a right Ore set. Then $(aS)T^{-1}$ is a right ideal for any $a\in S$.
\end{lemma}
 \begin{proof}
 Let $x\in (aS)T^{-1}$ and $b\in S$. We need to show that $xb\in (aS)T^{-1}$. By definition there exists $t\in T$ with $xt\in aS$.
 Since $T$ is a right Ore set, there exist $t'\in T$ and $b'\in S$ such that $bt'=tb'$. Therefore, $xbt'=xtb'\in aS$ and $xb\in (aS)T^{-1}$.
\end{proof}

In the following definition $P$ is a completely prime right ideal and $T=S-P$.

\begin{definition}\label{Def 3.2} A semigroup $S$ is \emph{right comparable} with respect to $P$ if for every $a,b\in S$
 one of the following three conditions holds: $aS\subseteq bS$, $bS\subseteq aS$ or $(aS)T^{-1}=(bS)T^{-1}$.
 We simply say that $S$ is a right $P$-comparable semigroup.
\end{definition}

\begin{example}\label{Ex 3.3}
\begin{itemize}
\item[(i)]
The set $S=\{0,e,f,ef,1,x,x^2,\dots\}$ with $ex=xe=x=fx=xf=efx$, $ef=fe$ and $e=e^2, f=f^2$ is a semigroup. 
The set $P=\{x^n: n\in \Nat\}\cup\{0\}$ is a completely prime ideal of $S$.
 Now consider the sets $(eS)=\{0,e,ef,x,x^2\dots\}$,  $(fS)$=$\{0,f,fe,x,x^2\dots\}$ and  $(efS)$=$\{0,ef,x,x^2\dots\}$. Let $a, b\in S$. 
If $a=x^i$ and $b=x^j$ then $aS\subseteq bS$ or $bS\subseteq aS$. 
Moreover, if $a=e$ and $b=f$ then $aS\nsubseteq bS$ and also $bS\nsubseteq aS$ but $(eS)T^{-1}=(eS)\cup\{f\}$ and $(fS)T^{-1}=(fS)\cup\{e\}$. 
This shows that $(fS)T^{-1}=(eS)T^{-1}$. In the case $a=ef$ and $b=e$ or $b=f$ we have $aS\subset bS$. Hence $S$ is a right $P$-comparable semigroup.

\item[(ii)]
In general, let $H$ be a right chain semigroup with identity $1$ and zero element $0\neq 1$. 
As in part (i), we can consider $S=H\cup \{e,f,ef\}$ with the relation $eh=fh=efh=h$ for all $h\in H-\{1\}$, and $e^2=e$, $f^2=f$, $ef=fe$. 
Then $J(H)$ is a completely prime ideal of $S$ and $S$ is a right $J(H)$-comparable semigroup. Furthermore, $S$ is not a right chain semigroup.
We do not know if for a right $P$-comparable semigroup $S$ there is a one sided chain semigroup $H$ with $P=J(H)$.
\end{itemize}
\end{example}

For the remainder of the paper completely prime ideals are always two sided ideals, unless state otherwise. 

Note that if $T\subseteq T'$ are multiplicatively closed subsets of $S$, 
then we have $(aS)T'^{-1}\subseteq (aS)T^{-1}$ for any $a\in S$.
We can therefore easily see that if $P'\subseteq P$ are completely prime 
ideals and if $S$ is $P$-comparable, then $S$ is a $P'$-comparable semigroup.

\begin{lemma}\label{Le 3.4}
Assume that $S$ is a right $P$-comparable semigroup, where $P$ is a completely prime ideal of $S$. 
Then $P$ is a right waist and the set of all completely prime ideals contained in $P$ are linearly ordered. 
Hence $N(S)$, the intersection of all completely prime ideals of $S$, is completely prime and a right waist.
\end{lemma}
\begin{proof}
Let $I\not\subseteq P$ be a right ideal of $S$ and suppose $a\in P$, 
$b\in I-P$. Since $bS\subseteq aS$ and 
$(aS)T^{-1}=(bS)T^{-1}$ imply a contradiction, we must have 
$aS\subseteq bS$ and so $P\subset I$. 
The rest follows from the remark preceding the lemma.
\end{proof}

\begin{proposition}\label{Pr 3.5}
Let $P$ be a proper completely prime ideal of $S$ and let $T=S-P$. Then the following conditions are equivalent:
\begin{itemize}

\item[(i)]
$S$ is a right $P$-comparable semigroup.

\item[(ii)]
For all $a,b\in S$, either $aS\subseteq bS$ or $(bS)T^{-1}\subseteq (aS)T^{-1}$.

\item[(iii)]
For all $a,b\in S$, either $aS\subseteq bS$ or $bS\subseteq (aS)T^{-1}$.

\item[(iv)]
$T$ is a right $Ore$ set and for all $a,b\in S$ we have $aS\subseteq bS$ or $b\in (aS)T^{-1}$.

\item[(v)]
For all $a\in S$, $(aS)T^{-1}$ is a right ideal and a right waist.
\end{itemize}
\end{proposition}
\begin{proof}
The proof is the same as that in ring theory, see \cite[Proposition 1.4]{FSA99}.
\end{proof}

Recall that a completely prime ideal $P$ of $S$ is a right waist if and only if $aP=P$ for every $a\in S-P$.

\begin{theorem}\label{Thm 3.6}
Let S be a right P-comparable semigroup. Then:
\begin{itemize}

\item[(i)]
 Any semiprime right ideal of $S$ contained in $P$ is a prime right ideal and a right waist.

\item[(ii)]
 The set of all prime right ideals contained in $P$ is linearly ordered by inclusion, and the prime radical $\beta(S)$ is a prime ideal and a right waist.

\item[(iii)]
 An ideal $Q$ of $S$ contained in $P$ is completely prime if and only if $Q$ is completely semiprime.
\end{itemize} 

\end{theorem}
\begin{proof}
To prove (i), let $Q$ be a semiprime right ideal contained in $P$ and $a,b\in S$. Suppose $aSb\subseteq Q$ and $a\notin Q$. If $a\in (bS)T^{-1}$,
 then there exists $t\in T, s\in S$ such that $at=bs$. Thus $atSat=atSbs\subseteq aSbs\subseteq Qs\subseteq Q$, and so $at\in Q$. 
Since $t\notin P$, $a\in P$, and $aSa\subseteq aP=atP\subseteq Q$, a contradiction.
 Consequently $a\notin (bS)T^{-1}$ and, by Proposition~\ref{Pr 3.5} part (iv), we have $bS\subseteq aS$. Hence $bSb\subseteq aSb\subseteq Q$ and $b\in Q$.

Let $b\in S$ and $bS\nsubseteq Q$. To prove that $Q$ is a right waist, it is enough to show that $aS\subseteq bS$ for every $a\in Q$. 
If $aS\nsubseteq bS$ for some $a\in Q$ then $b\in (aS)T^{-1}$ by part (iv) of the Proposition~\ref{Pr 3.5}.
 Therefore, $bt\in (aS)\subseteq Q\subseteq P$ for some $t\in T$. Since $t\notin P$, we have $b\in P$. 
Hence $bSb\subseteq bP=btP\subseteq Q$, so that $b\in Q$, a contradiction. 

The proof of (ii) follows from (i). To prove (iii), let $Q$ be completely semiprime and $ab\in Q$. 
Since $(bSa)^2\subseteq Q$ implies $(bSa)\subseteq Q$ and recalling Part (i) we have $b\in Q$ or $a\in Q$.
\end{proof}

\begin{lemma}\label{Le 3.7}
Let $S$ be a right $P$-comparable semigroup. If $I\subseteq P$ is a right ideal and a right waist, then $aI$ is also right waist for every $a\in S$.
 In particular, $aP$ is a right waist for every $a\in S$.
\end{lemma}

\begin{proof}
Assume that $x\notin aI$. Since $S$ is right $P$-comparable, we have $xS\subseteq aS$ or $a\in (xS)T^{-1}$. 
In the first case, there exists $s\in S-I$ with $x=as$, and so $I\subset sS$. Hence $aI\subseteq asS=xS$. 
In the second case, $at=xs$ for some $t\in T$ and $s\in S$. Therefore, $aI\subseteq aP=atP=xsP\subseteq xS$. 
\end{proof}
\begin{theorem}\label{Thm 3.8}
Let $S$ be right $P$-comparable with left cancellation law. Then $aP=bP$ if and only if $(aS)T^{-1}=(bS)T^{-1}$.
\end{theorem}

\begin{proof}
If $a$ and $b$ are not comparable, then $(aS)T^{-1}=(bS)T^{-1}$. 
We shall show that $aP=bP$. 
If $aP\neq bP$, we can assume   
$aP\subset bP$, because both $aP$ and $bP$ are right waists. Since $b\in (aS)T^{-1}$, there exist $t\in T$, 
$s\in S$ with $bt=as$. Hence $bP=btP=asP\subseteq aP$, a contradiction. 
Next we suppose that $a$ and $b$ are comparable. We can assume $aS\subseteq bS$. 
Therefore, $(aS)T^{-1}\subseteq (bS)T^{-1}$ and $aP\subseteq bP$. 
If $aP=bP$ and $(aS)T^{-1}\subset (bS)T^{-1}$ 
then $at\in bS$ for some $t\in T$ and $b\notin (aS)T^{-1}$. 
Thus $at=bp$ for some for some $p\in P$. 
Hence $aP=atP=bpP\subset bP$, a contradiction which shows that $a$ and $b$ can not be comparable. 
Conversely, if $(aS)T^{-1}=(bS)T^{-1}$ then, by the first part, we have $aP=bP$.
\end{proof}
We say that a semigroup $S$ is \textit{weak right $P$-comparable} if for every $a,b\in S$ one of the following three conditions holds:
 $aS\subseteq bS$, $bS\subseteq aS$ or $aP=bP$, where $P$ is completely prime. 
In ring theory, Ferrero and Sant'~Ana presented an example of right weak $P$-comparable 
ring that is not right $P$-comparable, see \cite[Example 2.6]{FSA99}. 
The above Theorem shows that in a left cancellation semigroup $S$, the notions of weak right $P$-comparable and right $P$-comparable coincide.
\begin{corollary}\label{Co 3.9}
If $S$ is a left cancellation semigroup, then $S$ is right $P$-comparable if and only if $S$ is right weak $P$-comparable. 

\end{corollary}
\begin{proof}
The proof is straight forward from the definition of right $P$-comparable and Theorem~\ref{Thm 3.8}.
\end{proof}

Let $S$ be a right $P$-comparable semigroup with left cancellation law. As in ring theory, we can define an equivalence relation on $S$.
 For elements $a, b\in S$ we put $a\sim b$ if and only if $aP=bP$. For $a\in S$, let $[a]=\cup_{b\sim a}bS$. 
The following proposition shows that $[a]$ is a right ideal and a right waist.

\begin{proposition}\label{Pr 3.10}
Let $S$ be a right $P$-comparable semigroup with left cancellation law. Then $[a]=(aS)T^{-1}$ is a right waist.
 Moreover, $(bS)T^{-1}=[a]$ for every $b\sim a$.
\end{proposition}

\begin{proof}
The proof by using Theorem~\ref{Thm 3.8} and  part (v) of Proposition~\ref{Pr 3.5} 
is the same as that in ring theory see \cite[Proposition 3.6]{FT93}.
\end{proof}
\begin{lemma}\label{Le 3.11}
Let $S$ be a right $P$-comparable semigroup with left cancellation law, and $I$ a right ideal with $P_r(I)=P\subseteq J=J(S)$. 
Then $I=\cup_{a\in I}(aS)T^{-1}$, where $T=S-P$, and $I$ is a right waist.
\end{lemma} 
\begin{proof}
Let $a\in I$, if $x\in (aS)T^{-1}$ then $xt\in aS$, for some $t\in T$. Hence $xt\in I$ and $t\not\in P_r(I)$.
 Therefore, $x\in I$ by definition of $P_r(I)$. Thus $(aS)T^{-1}\subseteq I$ for every $a\in I$ so that $\cup_{a\in I}(aS)T^{-1}\subseteq I$. 
The converse is clear. By Proposition~\ref{Pr 3.10} and the fact that the union of right waists is a right waist, $I$ is a right waist.
\end{proof}
\begin{corollary}\label{Co 3.12}
Let $I$ and $P$ be as the Lemma~\ref{Le 3.11}. Then $I=\cap_{a\notin I}aP=\cap_{a\notin I}aJ$, where $J=J(S)$ is the set of non units of $S$.
\end{corollary}
\begin{proof}
By Lemma~\ref{Le 3.11} $I$ is a right waist. By Corollary~\ref{Co 2.14} the remainder of the proof is obvious.
\end{proof}

The following theorem is an analogue of a characterization result for $D$-rings, due to Ferrero and Torner, see in 
\cite[Theorem 3.9]{FT93}. 
\begin{theorem}\label{Thm 3.13}
Suppose $S$ is a right $P_r(I)$-comparable semigroup with left cancellation law and $I$ is a nonzero right ideal of $S$. 
Then the following conditions are equivalent:
\begin{itemize}
\item[(i)]
 $P_r(I)\subseteq J=J(S)$.

\item[(ii)]
 There exists a completely prime right ideal $P$ contained in $J$ and a subset $V$ of $S$ such that $I=\cap_{a\in V}aP$.

\item[(iii)]
 There exist a completely prime right ideal $P$ contained in $J$ and a subset $V'$ of $S$ such that $I=\cup_{a\in V'}(aS)T^{-1}$, where $T=S-P$.
 \end{itemize}
 Furthermore, under the equivalent conditions above, $I$ is a right waist.
 
\end{theorem}
\begin{proof}
The equivalence of (i) and (ii) follows from Corollary~\ref{Co 3.12}. That (i) implies (iii) follows from Lemma~\ref{Le 3.11}.
 To see that (iii) implies (i), let $x\in P_r(I)$. Then there exists $s\in S-I$ with $sx\in I$, and $sx\in (aS)T^{-1}$ for some $a\in V'$.
 Hence $sxt\in aS$ for some $t\in T$. If $x\notin P$ then we have $s\in (aS)T^{-1}\subseteq I$, a contradiction. Thus $x\in P$ and $P_r(I)\subseteq P\subseteq J$. 
\end{proof}

\begin{lemma}\label{Le 3.14}
Let $S$ be a right $P$-comparable semigroup with left cancellation law. If $Q\subseteq P$ is a completely prime ideal, then $P_r(aQ)=Q$.
\end{lemma}

\begin{proof}
It is clear that $a\notin aQ$ and $ab\in aQ$ for every $b\in Q$. Hence $Q\subseteq P_r(aQ)$. Suppose that $x\in P_r(aQ)-Q$. 
Then there exists $c\in S-aQ$ with $cx\in aQ$. Thus $cx=aq$ for some $q\in Q$. Since $c\notin aQ\subseteq aP$ we have $(cS)T^{-1}\neq (aS)T^{-1}$.
 From Definition~\ref{Def 3.2}, we have $aS\subset cS$ or $cS\subseteq aS$. If $a=cr$ for some $r\in S$, then $aQ=crQ\subseteq cQ$. If $c=at$ for some $t\in S-Q$, then $cQ=atQ=aQ$.
 Therefore we have $aQ\subseteq cQ$ in all the cases. 
From $cx=aq$ and $aQ\subseteq cQ$, we have $cx=aq=cp$ for some $p\in Q$. Since $S$ is left cancellation and $cp=cx\neq 0$, we have $x=up$ for some unit $u\in S$, a contradiction. 
\end{proof}

The analogous for right cones of the following proposition was proved by Brungs and Torner in 
\cite[Proposition 1.11]{BT98}. 

\begin{proposition}\label{Pr 3.15}

Let $t\in P$ and $S$ a right $P$-comparable semigroup with left cancellation law. 
Then $Q=\cap t^nS$ is a prime right ideal and a right waist, if $t^nS\neq (0)$ for all $n\in \Nat$. 
In addition, if $Q$ is a two-sided ideal, then it is completely prime.
\end{proposition}

\begin{proof}
Since $t\in P$ and $S$ is left cancellation, we have $t^{n+1}S\subset t^nS$. Also $t^n\neq 0$ for all $n\geq 1$ by assumption.
 Hence $Q\subset t^nS$ for all $n$. If $x\notin Q$, then there exits an $n$ such that $x\notin t^nP$, and $t^nP\subset xS$, since $t^nP$ is a right waist. 
Thus $t^{n+1}S\subseteq t^nP\subseteq xS$. Hence we have $Q\subset t^{2n+2}S\subseteq t^{n+1}xS\subseteq xSxS$, and $xSx\nsubseteq Q$ follows. 
By part (i) of Theorem 3.6, $Q$ is prime right ideal and a right waist. 

Let $Q$ is a two-sided ideal and $x\notin Q$. Since $Q$ is a right waist we have $Q\subset xS$. 
Thus there exists an $n$ with $t^n=xa$ for some $a\in S$, and so $t^{2n}=xaxa$. We compare $ax$ and $x$. 
By (iv) of Proposition~\ref{Pr 3.5}, we have $axS\subseteq xS$ or $x\in (axS)T^{-1}$, where $T=S-P$. 
In the first case, we have $ax\in xS$, $t^{2n}=x^2b$ for some $b\in S$ and $x^2\notin Q$. By part (iii) of Theorem~\ref{Thm 3.6}, $Q$ is completely prime. 
In the second case there exist $t'\in T$ and $s\in S$ such that $xt'=axs$. If $xt'\in Q$, then $xP=xt'P\subseteq Q$ and $xSP\subseteq Q$. 
Hence $xS\subseteq Q$ or $P\subseteq Q$, a contradiction. This contradiction shows that $xt'\notin Q$. 
Since $Q$ is  a two-sided ideal, $xs\notin Q$, and so $xs\notin (t^mS)$ for some $m$. 
Again by (iv) of Proposition~\ref{Pr 3.5}, $t^m\in (xsS)T^{-1}$ and hence $t^mk=xsr$ for some $k\in T$ and $r\in S$.
 From $xt'=axs$ and $t^n=xa$, we have $x^2t'r=x(axs)r=(xa)(xsr)=t^nt^mk=t^{n+m}k$. Since $t^{n+m}k\notin Q$, we have $x^2\notin Q$.  

\end{proof}

\begin{example}\label{Ex 3.16}
\begin{itemize}

\item[(i)] Let $H$ be the semigroup $H=\{0\}\cup\{t^rx^n: 0\leq r\in {\Real},~n\in {\Z}\}$ with defining relation $xt^r=t^{2r}x$.
 If $I$ a nonzero two sided ideal of $H$, then $I=J(H)=\{t^rx^r:r>0\}$, since for any $t^rx^n\in I$ and positive integer $m$, we can write $t^rx^n=x^{m}t^{r/2m}x^{n-m}$.
 Let $Q=\cap t^nH$. It is clear that $Q$ is not a two sided ideal. Therefore, $Q$ is only a right prime ideal.

\item[(ii)] Proposition~\ref{Pr 3.15} is not true if $t\notin P$. 
To see this, let $S$ be as in Example~\ref{Ex 3.3}.  
Then $Q=\cap(ef)^nS=\{0,ef,x,x^2,\dots\}$ is a two sided ideal such 
that $ef\in Q$ but $e,f\notin Q$. 
Therefore $Q$ is not completely prime. 
\end{itemize}
\end{example}

\section{Prime segments of right $P$-comparable semigroups}\label{Sec 4}

Recall that a proper right ideal $A$ of a semigroup $S$ is called a comparizer (strongly comparizer) right ideal
 if for every $a,b\in S$, either $aS\subseteq bS$ or $bA\subseteq aS$ ($aS\subseteq bS$ or $bA\subseteq aA$).

Let $S$ be a right $P$-comparable semigroup. 
By Definition~\ref{Def 3.2}, Lemma~\ref{Le 3.4} 
and Theorem~\ref{Thm 3.8}, $P$ is a strongly comparizer ideal and right waist.
 Also, by Theorem~\ref{Thm 3.6}, any semiprime right ideal of $S$ contained in $P$ is a prime right ideal and a right waist.
 Therefore, any semiprime right ideal contained in $P$ is comparable with any ideal of a right $P$-comparable semigroup. 
It seems natural that right $P$-comparable semigroups with left cancellation law have many properties in common with right cones.
In this section we shall investigate those properties of ideal theory of cones that can be carried over right $P$-comparable semigroups. 
Throughout this section all semigroups have a left cancellation law, unless stated of otherwise. 
\begin{definition}\label{Def 4.1}
Let $S$ be a left cancellation semigroup and $P_2\subset P_1$ be completely prime ideals of $S$ such that there are no further completely prime ideals between $P_2$ and $P_1$.
Then we say that $P_2\subset P_1$ is a \emph{prime segment} of $S$. 
If $P_1$ is the minimal completely prime ideal of $S$, then $\varnothing \subset P_1$ is also considered a prime segment.
\end{definition}
\begin{definition}\label{Def 4.2}
A prime ideal $Q$ of a semigroup $S$ is called \emph{exceptional} if $Q$ is not completely prime.
\end{definition}
Let $A$ be an ideal of a semigroup $S$. An ideal $I$ (respectively, an element $s$) of $S$ is said to be \emph{$A$-nilpotent} if $I^n\subseteq A$ (respectively, $s^n\in A)$ for some $n\in \Nat$.

\begin{definition}\label{Def 4.3}
Let $A\subset B$ be ideals of a semigroup $S$ such that there are no further proper ideals between $A$ and $B$. Then we say that $B$ \emph{is minimal over} $A$.
\end{definition}
The following lemma is an extension of the pairing Lemma of \cite[Lemma 1.12]{BT98} to right $P$-comparable semigroup. 
\begin{lemma}\label{Le 4.4}(Pairing Lemma)
Let $S$ be a right $P$-comparable semigroup and $Q\subset P$ an exceptional prime ideal of $S$. Then there exists a unique waist ideal $D\supset Q$ that is minimal over $Q$. Furthermore, $D$ is an idempotent.
\end{lemma}
\begin{proof}
Set $D=\cap\{I: Q\subset I$ and $I$ is a right waist ideal of $S$\}. $D$ is nonempty since $Q\subset P$, and $D$ is a right waist. 
By part (i) of Lemma~\ref{Le 2.12} we have $Q\subset P_r(Q)$. From Lemma~\ref{Le 2.12} (ii), it follows that $Q\subset D$, and it is clear that $D$ is minimal over $Q$.
Since $Q$ is prime, $D$ is an idempotent.
\end{proof}

The following lemma is a generalization of \cite[Lemma 1.3]{BT98} and \cite[Lemma 17]{FMSA05} to right $P$-comparable semigroups. 
\begin{lemma}\label{Le 4.5}
Let $Q$ be an exceptional prime ideal of a right $P$-comparable semigroup $S$, and let $D$ be the idempotent right waist of $S$, minimal over $Q$. Then there exists an element $a\in D-Q$ such that 
 $Q\subset \cap_{n\in\Nat}a^nS$. In particular, there exists a non $Q$-nilpotent element in $D-Q$.
\end{lemma}
\begin{proof}
Since $Q$ is an exceptional prime ideal of $S$, there 
exists $b\in S-Q$ with $b^2\in Q$ by part (iii) of the Theorem~\ref{Thm 3.6}. 
Set $W=\{c\in S:~\cap_{n\in \Nat}c^nS\subseteq Q\}$. 
Suppose that $c\in W-D$. Then, by part (i) of Lemma~\ref{Le 2.2}, we have $D=cD$. Thus $D=c^nD\subseteq c^nS$ for all $n$. Hence $Q\subset D\subseteq \cap_{n\in \Nat}c^nS$, a contradiction. 
Thus $W\subseteq D$. 
If $b\in cbS$ for some $c\in W$, then $b=cbs=c^nbs^n$ for some $s\in S$ and all $n$. Therefore, $b\in c^nbs^nS\subseteq c^nS$ for all $n$, and so
$b\in \cap_{n\in \Nat}c^nS\subseteq Q$, a contradiction. 
By (iv) of Proposition~\ref{Pr 3.5}, $cb\in (bS)T^{-1}$ for all $c\in W$. 
Hence $Wb\subseteq (bS)T^{-1}$. If $W=D$, then 
$bDb=bWb\subseteq b(bS)T^{-1}\subseteq Q$, and so $b\in Q$, a contradiction.
Hence $W\subset D$. Let $d\in D-W$. Then 
$\cap_{n\in \Nat}a^nS\nsubseteq Q$, and since 
$Q$ is a right waist, we have $Q\subset \cap_{n\in \Nat}a^nS$. 
\end{proof}

See \cite[Lemma 2.4]{BS95} for a similar result of the following lemma.
\begin{lemma}\label{Le 4.6}
Let $S$ be a right $P$-comparable semigroup and
\[
\alpha (P)=\{P'\subset P:~P' \text{ is semiprime}\}.
\]
Then: 
\begin{itemize}

\item[(i)]
 The set $\alpha (P)$ is totally ordered.

\item[(ii)]
 The set $\alpha (P)$ is closed under union and intersection. 

\item[(iii)]
 The set $\alpha (P)$ has a greatest lower bound $P_0\in \alpha (P)$.  

\item[(iv)]
 For any completely semiprime ideal $P'$ with $P'\subset P$, there exists a completely prime ideal $P_0$ of $S$
 such that $P'\subseteq P_0$ and $P_0\subset P$ is a prime segment.
 \end{itemize}
\end{lemma}
\begin{proof}
 All of (i), (ii) and (iii) follow directly from Theorem~\ref{Thm 3.6}. To prove (iv) let $a\in P-P'$ and set $L(a)=\{P_i\subset P : a\notin P_i$ and $P_i$ is completely semiprime\}.
 By Theorem~\ref{Thm 3.6}, every $P_i$ in $L(a)$ is completely prime and the set $L(a)$ is linearly ordered by inclusion. 
Therefore, the ideal $P_0=\cup_{P_i\in L(a)} P_i$ has the desired properties. 
\end{proof}
\begin{definition}\label{Def 4.7}
A prime segment $P_2\subset P_1$ of a semigroup $S$ is called \emph{Archimedean}, if for every $a\in P_1-P_2$
there exists an ideal $I\subseteq P_1$ with $a\in I$ and $P_2=\cap I^n$. It is called \emph{simple} 
if there are no further two-sided ideals of $S$ between $P_2$ and $P_1$. It is called \emph{exceptional} 
if there exists a prime ideal $Q$ of $S$ with $P_2\subset Q\subset P_1$ and no further two-sided ideal exists between $P_1$ and $Q$. 
\end{definition}

In the following we give a characterization theorem for prime segments, generalizing a similar of right cones, see \cite[Theorem 1.14]{BT98}.
\begin{theorem}\label{Thm 4.8}
Let $S$ be a right $P_1$-comparable semigroup, and let $P_2\subset P_1$ be a prime segment of $S$. Then exactly one of the following alternatives occurs:
\begin{itemize}

\item[(i)]
The prime segment $P_2\subset P_1$ is Archimedean.

\item[(ii)]
The prime segment $P_2\subset P_1$ is simple, i.e. there are no further ideals between $P_1$ and $P_2$.

\item[(iii)]
There exists a prime ideal $Q$ with $P_2\subset Q\subset P_1$ and no further ideal between $P_1$ and $Q$ exist.
Moreover, $P_2=\cap_{n\in \Nat}Q^n$. 
\end{itemize}
\end{theorem}

\begin{proof}
Let $L(P_1)=\cup I$, the union of ideals $I$ of $S$ properly contained in $P_1$. If $L(P_1)=P_2$,
then the prime segment is simple. 

Assume that $P_2\subset L(P_1)\subset P_1$ and $P_1=P_1^2$. 
Let $A$ and $B$ be ideals of $S$ with $L(P_1)\subset A, B$. Since $P_1$ is a waist, and $A, B\nsubseteq P_1$, we have $P_1\subseteq A$ and $P_1\subseteq B$. 
Then $P_1=P_1^2\subseteq AB$, and $L(P_1)$ is exceptional. Let $L(P_1)=Q$. Since $P_2$ is prime and a waist, 
we have $P_2\subset Q^n$ for any $n$, which shows that $Q$ is nonnilpotent. 
By part (iv) of Theorem~\ref{Thm 2.4}, $\cap_{n\in \Nat}Q^n$ is a completely prime ideal, and so $P_2=\cap_{n\in \Nat}Q^n$.

Assume that $P_1\neq P_1^2$ or $L(P_1)=P_1$, then $\cap P_1^n=P_2$, or for any $a\in P_1-P_2$ there exists an ideal $I$ with $a\in I$ and $P_2\subset I\subseteq P_1$. 
Since $I$ is nonnilpotent, again using part (iv) of Theorem~\ref{Thm 2.4}, we have $\cap_{n\in \Nat}I^n=P_2$. Therefore  the prime segment is Archimedean. 
\end{proof}

\begin{example}\label{Ex 4.9}
\begin{itemize}
\item[(i)]
 Let $S=\{0,1,x_1,x_2,x_3,\dots\}$ be a semigroup with the relation
$x_ix_j=x_{\min(i,j)}$. Then for each $i$ the ideal generated by $x_i$ is 
a completely prime ideal and also an idempotent. For each $i$,   
$Sx_iS\subset Sx_{i+1}S$ is a simple prime segment.

\item[(ii)]
Let $=\{0,1,e,f,ef,x,x^2,\dots\}$ and $P=\{0,x^n: n\in \Nat\}$ as in Example~\ref{Ex 3.3}.
Then $(0)\subset P$ is Archimedean prime segment.

\item[(iii)] 
Let $S=\{0,1,x_1,x_2,x_3,\dots,x_n\}$ is a semigroup with relation $x_ix_j=\delta _{i,j}x_j$, where $\delta _{i,j}$ is the Kronecker delta.
Then for each $1\leq i\leq n$, the set $P_i=S-\{1,x_i\}$ is a minimal completely prime ideal. Therefore, $\varnothing \subset P_i$ is a prime segment. 
This prime segment is neither simple nor Arthimedean and not exceptional. 
In this semigroup all ideals are idempotents but not necessarily completely prime. 
This example shows that the cancellation law is necessary for Theorem~\ref{Thm 4.8} and also for part (v) of the Theorem~\ref{Thm 2.4} to hold.

\item[(iv)]
 Let $H=\{t^rx^n:~r\geq 0, n\in {\Z}\}\cup\{0\}$. 
As mentioned in Example~\ref{Ex 3.16} (i) the Jacobson radical 
$J(H)=\{t^rx^n: r>0\}$ and $\{0\}$ are the only two-sided ideals of $H$. 
This kind of semigroup is called \emph{nearly simple}. Let 
\[
U=\left\{\begin{pmatrix} a&b\\0&a^{-1}\\
\end{pmatrix}:~b,0\leq a\in {\Real}\right\}.
\]
We take for $H'$ the subsemigroup $\{t^ru: 0\leq r,~u\in U\}$ of the group $SL(2, \Real)$.
Brungs and Dubrovin showed in \cite{BD03} that  $t^{\pi}H'$ is a prime ideal that is not completely prime. 
Therefore, the prime segment $(0)\subset J(H')$ is exceptional.  
\end{itemize}
\end{example}

Recall that a prime segment $P_2\subset P_1$ is called locally (right) invariant if ($P_1a\subseteq aP_1$) $P_1a=aP_1$ holds for every $a\in P_1-P_2$.

The following Lemma is an extension of \cite[Corollary 1.15]{BT98} to $P$-comparable semigroups.
\begin{lemma}\label{Le 4.10}
If $S$ is a right $P_1$-comparable semigroup with left cancellation law, then every locally invariant prime segment $P_2\subset P_1$ is Archimedean. 
\end{lemma}
\begin{proof}

By using Proposition 3.15 and the fact that all prime ideals contained in $P_1$ also are right waists, the proof is the same as for right cones
see \cite[Corollary 1.15]{BT98}. 
\end{proof}
The converse of the Lemma~\ref{Le 4.10} also appears to be true, but we have not been able to prove this.

\bibliographystyle{amsplain}

\end{document}